
\input amstex.tex 
\documentstyle{amsppt} 
\NoRunningHeads 
\baselineskip=12pt 
\parskip=5pt 
\parindent=18pt 
\topskip=10pt 
\leftskip=0pt 
\pagewidth{30pc} 
\pageheight{45pc} 
\topmatter 
\title 
{\bf On Endomorphisms of Algebraic Surfaces} 
\endtitle 
\author D. -Q. Zhang 
\endauthor 
\address 
\newline \noindent 
Department of Mathematics, National University of Singapore, 
2 Science Drive 2, 
\newline \noindent 
Singapore 117543 
\endaddress 
\par \noindent 
\email 
\newline \noindent 
matzdq$\@$math.nus.edu.sg 
\endemail 
\dedicatory 
\enddedicatory 
\abstract 
In these notes, we consider self-maps of degree $\ge 2$ on a weak del Pezzo
surface $X$ of degree $\le 7$. We show that there are exactly $12$ such
$X$, modulo isomorphism. In particular, $K_X^2 \ge 3$, and if $X$ has one
self-map of degree $\ge 2$ then for every positive integer $d$
there is a self-map of degree $d^2$ on $X$.
\endabstract 
\endtopmatter 

\document 
\head Introduction 
\endhead 

We work over the complex numbers field ${\bold C}$. 
By a {\it self-map} $f$, we mean an algebraic morphism from an algebraic variety $X$
onto itself. Invertible self-maps (e.g. automorphisms and hence of degree $1$)
have been intensively studied by mathematicians. However, those $f$ of
degree $\ge 2$ are not known much yet.

\par
We remark that the study of self-maps on singular varieties can not be
reduced to that on smooth varieties. Indeed, the multiplication by an
integer
$n \ge 2$ on an abelian surface $A$ induces a self-map $f$ of degree $n^4$ on
the Kummer
surface $K$ ($= A$ modulo the involution $\text{\rm -id}$). However,
the induced self {\it rational} map $f'$ on the minimal resolution of the Kummer surface
$K$ is not a morphism for the $f$-pre image of the 16 singularities on $K$
has more than these $16$ points.

\par
Self-maps of open algebraic varieties are actually more difficult. For example,
the Jacobian conjecture, which has been open for 50 more years and is still
open, asserts that every etale endmorphism (not necessarily surjective)
of the affine $n$-space is an automorphism.
In these notes we will consider self-maps of projective varieties only.

\par
We mention some known results. A self-map of a variety of Kodaira
dimension
$\ge 0$ is unramified (Sommese's result).
In particular, there are no self-maps of degree $\ge 2$
on a surface with non-negative Kodaira dimension and
non-zero $c_1^2$ or $c_2$.
Beauville [B] proved that there are no self-maps of degree $\ge 2$ on
a smooth hypersurface of degree $\ge 3$ in a projective space, and that
if a del Pezzo surface $X$ has a self-map 
of degree $\ge 2$, then $K_X^2 \ge 6$.

\par
In view of the above results, one needs to study the self-maps on
algebraic varieties of Kodaira dimension $-\infty$. 
The first thing is naturally the study of those on Fano or even weak Fano varieties.
Here we say a smooth projective variety is {\it Fano} (resp. {\it weak Fano})
if its anti-canonical divisor is ample (resp. nef and big). 
A Fano (resp. weak Fano surface) is also called
a {\it del Pezzo} (resp. {\it weak del Pezzo}) surface.
A weak del Pezzo surface $X$ satisfies $1 \le K_X^2 \le 9$. 
Also, $K_X^2 = 9$ if and only if $X = {\bold P}^2$,  and
$K_X^2 = 8$ if and only if $X$ is a Hirzebruch surface 
$\Sigma_s$ of degree $s$ for some $0 \le s \le 2$. 

\par
After the preliminary version of these notes was written,
I was kindly informed by Professor Dolgachev that a recent thesis
of M. Dabija [Da] has shown that if a blow up of ${\bold P}^2$ has a self-map
of degree $\ge 2$, then it is the blow up of at most 3 non-collinear points
(and their infinitely near points), which is, in the case of weak del Pezzo
surfaces, part of our Theorem 2 (see also Lemma 1.1).
Actually, we prove that there are exactly $12$ isomorphism classes
of weak del Pezzo surfaces of degree $\le 7$, each of which has at least one
(and hence infinitely many) self-maps of degree $\ge 2$.
Now we state our main results.

\par \vskip 1pc 
{\bf Theorem 1.} {\it Let $X$ be a weak del Pezzo surface with 
$K_X^2 \le 7$ and a self-map $f$ which is not an automorphism. Then we have:}
 
\par \noindent 
(1) {\it We have $K_X^2 \ge 3$. Conversely, for every $3 \le d_x \le 7$ 
and every $1 \le d_f$, there is a weak del Pezzo surface $X$ with 
$K_X^2 = d_x$ and a self-map $f$ of degree $d_f^2$.} 

\par \noindent 
(2) {\it There is an integer $m > 0$ such that
$g := f^m$ is again a self-map of $X$ 
of degree $d^2$ and 
induced from a self-map of ${\bold P}^2$, also denoted by $g$, i.e., 
there is a birational morphism $\mu : X \rightarrow {\bold P}^2$ 
such that $g \mu = \mu g$. Moreover, $g^*E = dE$ for every negative 
curve $E$ on $X$ (see Theorems $2$ and $2.1$ for the description of such
$g$).} 

\par \vskip 2pc 
{\bf Theorem 2.} {\it Let $X$ be a weak del Pezzo surface 
with $K_X^2 \le 7$ and a self-map $f$ (which is not an automorphism) 
such that $\text{\rm Supp}(f^{-1}(E)) = E$ for every negative curve $E$ on $X$. 
Then we have:}

\par \noindent
(1) {\it We have $deg(f) = d^2$ for some integer $d \ge 2$.
There are a birational morphism $\mu : X \rightarrow {\bold P}^2$ 
and a self-map of ${\bold P}^2$, also denoted by $f$, such that 
$\mu f = f \mu$.}

\par \noindent 
(2) {\it There is a triangle of three non-concurrent lines $\sum_i L_i$ 
such that $\mu$ is the blow-up of some of the three intersection points 
of the triangle and their infinitely near points satisfying 
the conditions $(2a)-(2c)$ in Theorem $2.1$ in Section $2$, 
especially the reduced divisor $LP = \text{\rm Supp}(\mu^{-1}(\sum L_i))$ is a simple
loop of ${\bold P}^1$'s and all 
negative curves $E_i$ of $X$ are contained in $LP$.} 

\par \noindent 
(3) {\it The ramification divisor of $f$ is determined below, where 
$Z /(d-1)+ \sum E_i \sim_{\bold Q} LP$} 
$$K_X = f^*K_X + (d-1)\sum_i E_i + Z.$$ 

\par \noindent 
(4) {\it The dual graph of the loop $LP$ in $(2)$ determines uniquely the
isomorphism class 
of $X$. There are exactly twelve such loops (and hence exactly 
twelve isomorphism classes of surfaces $X$, i.e., $X(j)$ ($1 \le j \le 12$)
in Lemma $2.4$) satisfying the conditions $(2a)-(2c)$ 
and $K_X^2 \le 7$ (see Figures $1-12$); see Theorem $2.1$ below For more information on $\mu$ and
the converse part.} 

\par \vskip 1pc
{\bf Theorem 3.} 

\par \noindent
(1) {\it Let $X$ be a weak del Pezzo surface with a self-map
of degree $\ge 2$ and with $K_X^2 \le 7$.
Then $X$  is isomorphic to one of those $12$ isomorphism
classes $X(j)$ ($1 \le j \le 12$) given in Lemma $2.4$;  in particular,
$K_X^2 \ge  3$.}

\par \noindent
(2) {\it For every $1 \le j \le 12$ and every integer $d \ge 1$, there is a
self-map of degree $d^2$ on the surface $X(j)$.}

\par \vskip 1pc
{\bf Remark 4.}

\par \noindent
(1) Not every self-map of a surface has degree $d^2$ for some integer $d$.
Indeed, if $f : {\bold P}^1 \rightarrow {\bold P}^1$ is
the self-map given by $X_i \mapsto X_i^d$, then 
the self-map $f \times \text{\rm id}_C$ of ${\bold P}^1 \times C$,
where $C$ is a curve, has degree $d$.

\par \noindent
(2) Not every automorphism of a blown up surface of ${\bold P}^2$
is lifted from a projective transformation of ${\bold P}^2$;
see [ZD, Theorem 1].

\par \vskip 1pc
{\bf Question 5.}

\par \noindent
(1) Suppose that $f$ is a self-map of degree $\ge 2$ on a rational surface.
Find the condition for $f$ to be lifted from a self-map on ${\bold P}^2$.

\par \noindent
(2) Suppose that $f$ is a self-map of $X$ such that $f^{-1}(P) = \{P\}$ for some
point $P$. Find a necessary and sufficient condition on $f$ so that one can lift $f$ to a self-map on the blow-up of $X$ at $P$. A sufficient condition for
lifting is given in Lemma 1.2.

\head
Terminology and Notation
\endhead

\par \noindent
1. All surfaces (and also varieties) in the paper are assumed to be
projective unless otherwise specified.

\par \noindent 
2. A {\it self-map} $f$ of an algebraic variety $X$ is an algebraic morphism
from 
$X$ onto itself. 

\par \noindent 
3. A $(-n)$-{\it curve} on a smooth surface is a smooth rational curve of
self 
intersection $-n$. 

\par \noindent 
4. A curve $C$ on a surface is {\it negative} if $C^2 < 0$. 
A divisor on a surface is {\it negative} if its irreducible components have
negative intersection matrix. 

\par \noindent 
5. For a divisor $D$, we use $\#D$ to denote the number of irreducible
components in $\text{\rm Supp} D$. 

\head
Acknowledgement
\endhead

The author is very grateful to Professor Dolgachev for his encouragement
and valuable suggestions.
This work was started when the author visited 
Korea Institute for Advanced Study in December 2001. The author would like 
to thank the Institute, especially Professor J. Keum for the valuable
discussion.
This work was also partially supported by an Academic Reserach Fund of
National University of Singapore. 

\par 
After the paper was presented in the conference on Algebraic Geometry
on the occasion of Professor S. Iitaka's 60th birthday, Professors
E. Sato and K. Cho kindly pointed out that Professor N. Nakayama has
classified surfaces with a self map of degree $\ge 2$, in particular:
a rational surface has a self map of degree $\ge 2$
if and only if it is a toric surface. This agrees with our Theorem 2,
where the surface is a blow-up of points of the intersection of the
triangle and their infinitely near point so that the inverse of the 
triangle is a simple loop and hence the resulted surface is a toric surface.
The methods of Professor Nakayama's and ours are different.
According to Professor Nakayama's preprint, there are also works by Segami
(resp. Fujimoto-Sato)
on self-maps for irrational surfaces (resp. threefolds).

\head
Section 1. Examples and preparations
\endhead

\par \vskip 1pc 
{\bf Lemma 1.1.} {\it Let $X$ be a smooth surface with a self-map $f$ 
of degree $\ge 2$ so that every negative divisor on $X$ is contractible to a ${\bold
Q}$-factorial 
singularity (this is true if $X$ is weak del Pezzo; see Lemma $1.4$)
Then we have:} 

\par \noindent 
(1) {\it $f$ is a finite morphism.} 

\par \noindent 
(2) {\it Suppose that $E$ is a negative curve on $X$. Then $F =
\text{\rm Supp}(f^{-1}(E))$ is 
irreducible and negative.} 

\par \noindent 
(3) {\it Suppose that there are (at least one and) only finitely many
negative curves on $X$ (this is true if $X$ is weak del Pezzo by Lemma
$1.4$). Then 
there is a positive integer $m$ such that $g = f^m$ is a self-map of 
$X$ satisfying $deg(g) = d^2$ for some integer $d > 0$ 
and $g^*E = dE$ for every negative curve $E$ on $X$.} 

\par \vskip 1pc 
{\it Proof.} The proof of (1) is similar to that of (2). 
For (2), let $\mu : X \rightarrow Y$ be the algebraic contraction of $E$ 
to a point $P$ on the normal surface $Y$. 
The Stein factorization of $\mu f$ will decompose as $\mu f = g \nu$, 
where $g : Y' \rightarrow Y$ is a finite morphism. Now $\nu$ maps $F$ to
points lying over the point $P = \mu f(E)$. 
Hence $F$ is contractible (and hence negative; see [M]) and the Picard
number $\rho(Y') = \rho(X) - 
\# F$. On the other hand, $Y'$ dominates $Y$ and hence $\rho(Y') \ge \rho(Y)
= \rho(X) - 1$. 
Thus $F$ is irreducible and negative. 

\par 
For (3), note that by (2), $f^{-1}$ induces a bijection on the finite set
$\Sigma$ of all negative curves on $X$. 
Hence, for some $m > 0$,  the inverse of $g = f^m$ acts as identity on the
set $\Sigma$, i.e.,  $\text{\rm Supp}(g^{-1}(E)) = E$ for every negative curve $E$ on
$X$. Now (3) follows from the proof of Lemma 1.2. This proves the lemma.

\par \vskip 1pc 
{\bf Lemma 1.2.} 
\par \noindent 
(1) {\it Suppose that $f$ is a self-map on a smooth surface $X$. 
If $E$ is a $(-1)$-curve on $X$ such that $\text{\rm Supp}(f^{-1}(E)) = E$ 
and $\mu : X \rightarrow Y$ is the blow down of $E$ 
to a point $P$, then $Y$ has an induced self-map, also denoted by $f$ such
that 
$f \mu = \mu f$ and $f^{-1}(P) = \{P\}$. Moreover, $deg(f) = d^2$ and $f^*E
= dE$ 
for some positive integer $d$.} 

\par \noindent 
(2) {\it Suppose that $f$ is a self-map on a smooth surface $Y$. 
If $P$ is a point on $Y$ 
(so that two smooth curves $C_i$ intersect transversally at $P$ 
and $f^*C_i = dC_i$; this technic condition could probably be weakened) such that
$f^{-1}(P) = \{P\}$ 
and $\mu: X \rightarrow Y$ is the blow up of $P$ with $E$ the exceptional
divisor, 
then there is an induced self-map of $X$, also denoted by $f$, such that 
$\mu f = f \mu$ and $f^*E = dE$, where $deg(f) = d^2$.} 

\par \vskip 1pc 
{\it Proof.} The first part of (1) follows from the proof of Lemma 1.1 (2). 
For the second part of (1), let $f^*E = d'E$ and $f_*E = d''E$. Then $deg(f)
E = f_*f^*E = d'd''E$. 
Also $deg(f) E^2 = (f^*(E))^2 = (d')^2 E^2$. Thus $d' d'' = deg(f) = (d')^2$
and $deg(f) = d^2$ with $d = d' = d''$. 

\par 
For (2), suppose that $C_i = \{x_i = 0\}$ with $x_1, x_2$ the local
coordinates at $P$. 
Then the blow up of $P$ is locally $(x_i, x_j/x_i) \mapsto (x_i, x_j)$ with
$\{i, j\} = \{1, 2\}$, while $f$ is given by $f^* : x_i \mapsto x_i^d \times$
(a unit). 
Now it is clear that $f$ is liftible. This proves the lemma.

\par \vskip 1pc 
{\bf Lemma 1.3.} 

\par \noindent
(1) {\it A weak del Pezzo surface is a rational surface.} 

\par \noindent 
(2) {\it Let $X$ be a smooth rational surface and let $D$ be an effective
divisor on $X$. 
Suppose that $|K_X+D| = \emptyset$. Then $D$ is of simple normal crossing
and supported by smooth rational curves.} 

\par \noindent
(3) {\it Let $X$ be a smooth rational surface and let $0 \ne D$ be a reduced
divisor
in $|-K_X|$. Then either $D$ is irreducible and of arithmetic genus $1$,
or $D$ is reducible and every component of $D$ is isomorphic to ${\bold P}^1$.
Moreover, if $D$ is reducible and let $D = \sum D_j$ with $D_j$ irreducible,
then either $\#D = 2$ and $D_1$ and $D_2$ are tangent to each other (of
order $2$) at a point,
or $\#D = 3$ and $D_j$'s pass through a common point, or $D$ is a simple
loop.}

\par \vskip 1pc 
{\it Proof.} Note that the contraction of all $(-2)$-curves on $X$ 
will give rise to a Gorenstein del Pezzo surface (i.e., the anti-canonical 
divisor is ample and all singularities are Du Val), which is rational so
does $X$, see e.g. [GZ1, Lemma 1.3]. 
(2) is a consequence of the Riemann-Roch theorem applied to all sub-divisors
of $D$. 
For (3), let $B$ be the minimal sub-divisor of $D$, which is not a normal crossing
rational tree.
Then the Riemann-Roch theorem implies that $|K_X+B| \ne \emptyset$.
Hence $B = D \sim -K_X$. Also $D \sim -K_X$ implies that one of the cases in
(3) occurs. We need to say that there is such $B$.
Otherwise, $D$ is a disjoint union of $e$ rational trees.
So $0 = D . (D +K_X) = -2e$, a contradiction.
This proves the lemma.

\par \vskip 1pc 
{\bf Lemma 1.4.} {\it Let $X$ be a weak del Pezzo surface. 
Then we have:} 

\par \noindent 
(1) {\it If $C$ is a curve on $X$ with $C^2 < 0$, then $C$ is a $(-1)$ or
$(-2)$-curve.} 

\par \noindent 
(2) {\it There are only finitely many curves $C$ on $X$ with $C^2 < 0$.} 

\par \vskip 1pc 
{\it Proof.} For (1), since $-K_X$ is nef, we have $K_X . C \le 0$. 
Now (1) follows from the genus formula for $C$. 

\par 
(2) is well known. We prove it here for readers' convenience. 
Note that all $(-2)$-curves are perpendicular to the nef and big divisor 
$-K_X$, whence contractible altogether to Du Val singularities
by using the Hodge index theorem. So the number of all $(-2)$-curves on $X$ 
is bounded from above by $\rho(X) - 1$. 

\par 
For $(-1)$-curves $E_i$, note that $E_i' := E_i + (1/d)K_X$ is perpendicular
to $K_X$ 
where $d = K_X^2$ and hence can be written as 
$$E_i' = \sum_j a_j B_j$$ 
where $K_X, B_j$ form an orthogonal basis of $(Pic X) \otimes {\bold Q}$ 
(so $B_j^2 < 0$, with $B_j$ chosen to be integral). Since $1 + 1/d = -(E_i')^2 = \sum
a_j^2 (-B_j^2)$, 
the values $|a_j|$ 
are bounded from above. Since 
$a_k B_k^2 = E_i' . B_k \in
(1/d){\bold Z}$, 
the denominators of $a_j$ are bounded from above. Thus there are only
finitely 
many choices of $a_i$ and hence only finitely many cohomology classes
$[E_i']$. 

\par 
Suppose the contrary that there are infinitely many $(-1)$-curves $E_i$
on $X$. Then $[E_a'] = [E_b']$ for some $a \ne b$. So
$-1 - 1/d = (E_a')^2 = (E_a'  . E_b') = E_a . E_b - 1/d$
and $E_a . E_b = -1$, a contradiction.  So the lemma is true.

\par \vskip 1pc 
{\bf Example 1.5.} Let $f = \varphi_d : {\bold P}^2 \rightarrow {\bold P}^2$
be the morphism 
of degree $d^2$ ($d \ge 2$) given by 
$$[X:Y:Z] \rightarrow [X':Y':Z'], \,\,\, X' = X^d, \,\, Y' = Y^d, \,\, Z' = Z^d.$$
Let $L_i$ (resp. $L_i'$) be the lines defined respectively by $X = 0$, 
$Y = 0$ and $Z = 0$ (resp. $X' = 0$, $Y' = 0$ and $Z' = 0$). 
Then $f^*(L_i') = d L_i$ and 
$$K_{{\bold P}^2} = f^*(K_{{\bold P}^2}) + (d-1) \sum_i L_i.$$  
If we use $[U:V:W]$ as coordinates for both the domain and range of $f$ 
so that $f$ is a self-map and use $L_i$ to denote the lines defined by 
coordinates, then we have $f^*L_i = dL_i$ and $f^{-1}(P_i) = P_i$ 
where $P_1 = [1:0:0], P_2 = [0:1:0], P_3 = [0:0:1]$. 
Set $X_0 = {\bold P}^2$ and $LP_0 = \sum L_i$. 
Let $\mu_1 : X_1 \rightarrow {\bold P}^2$ be the blow up 
of one of $P_i$. 
By Lemma 1.2, we have an induced self-map $f$ on $X_1$ such that 
$f^*C = dC$ for every curve in the reduced divisor $LP_1 = \mu_1^{-1}(\sum_i L_i)$ 
which is a simple loop so that $K_{X_1} + LP_1 \sim 0$. Moreover, 
$$K_{X_1} = f^*K_{X_1} + (d-1)(LP_1).$$ 
Similarly, we let $\mu_2 : X_2 \rightarrow X_1$ be the blow up 
of an intersection point of the loop $LP_1$ and set the reduced divisor $LP_2 =
\mu_2^{-1}(LP_1)$. 
Thus we obtain $\mu_i : X_i \rightarrow X_{i-1}$ and a self-map $f$ on all
$X_i$ commuting 
with all $\mu_j$ such that 
$$K_{X_i} = f^*K_{X_i} + (d-1)(LP_i),$$ 
where the reduced divisor $LP_i$ is the inverse of the triangle $LP_0 = \sum
L_i$ 
and is again a simple loop so that 
$$K_{X_i} + LP_i \sim 0.$$

\head
Section 2. The proof of Theorems 1 and 2
\endhead

\par \vskip 1pc 
We will prove Theorem 2.1 below which includes Theorem 2 in the
Introduction. 

\par \vskip 1pc 
{\bf Theorem 2.1.} {\it Let $X$ be a weak del Pezzo surface 
with $K_X^2 \le 7$ and a self-map $f$ (which is not an automorphism) 
such that $\text{\rm Supp}(f^{-1}(E)) = E$ for every negative curve $E$ on $X$. 
Then we have:}

\par \noindent
(1) {\it We have $deg(f) = d^2$ for some integer $d \ge 2$. 
There are a birational morphism $\mu : X \rightarrow {\bold P}^2$ 
and a self-map of ${\bold P}^2$, also denoted by $f$, such that 
$\mu f = f \mu$.}

\par \noindent 
(2) {\it There is a triangle of three non-concurrent lines $\sum_i L_i$ 
such that $\mu$ is the blow-up of some of the three intersection points 
of the triangle and their infinitely near points satisfying} 

\par \noindent 
(2a) {\it The reduced divisor $LP = \mu^{-1}(\sum L_i)$ is a simple loop  of ${\bold P}^1$'s;
write $LP = \sum_i E_i + \sum_j P_j$ with $E_i^2 < 0$ and $P_j^2 \ge 0$,} 

\par \noindent 
(2b) {\it every negative curve on $X$ is contained in the loop $LP$, 
i.e., equal to one of $E_i$,} 

\par \noindent 
(2c) {\it every curve in the loop $LP$ has self intersection in
$\{-2,-1,0,1\}$.} 

\par \noindent 
(3) {\it We have $K_X+LP \sim 0$ and $f^*E_i = d E_i$;
the ramification divisor of $f$ is determined below} 
$$K_X = f^*K_X + (d-1)\sum_i E_i + Z$$ 
{\it where the effective integral divisor $Z$ is supported by curves 
of non-negative self intersection, such that $Z \sim (d-1) \sum
P_j$.}
 
\par \noindent 
(4) {\it The dual graph of the loop $LP$ in $(2)$ determines uniquely the
isomorphism class 
of $X$. There are exactly twelve such loops (and hence exactly 
twelve isomorphism classes of surfaces $X$,
i.e., $X_j$ ($1 \le j \le 12$) in Lemma $2.4$) satisfying the conditions
$(2a)-(2c)$ 
and $K_X^2 \le 7$ (see Figures $1-12$).} 

\par \noindent
(5) {\it Conversely, if $\nu : Y \rightarrow {\bold P}^2$ is the blow-up of 
some intersection points of a triangle of three non-concurrent lines 
$\sum L_i$ satisfying $(2a)-(2c)$ and $K_Y^2 \le 7$, then $Y$ 
is a weak del Pezzo surface with $K_Y^2 \ge 3$.
For every $d > 1$, this $Y$ has a self-map $g$ 
of degree $d^2$ (induced from a self-map $g$ of ${\bold P}^2$ 
so that $g \nu = \nu g$) such that $g^*C = dC$ for every curve in $LP$, where
the reduced divisor $LP$ is the simple loop $\nu^{-1}(\sum L_i)$, and  that} 
$$K_Y = g^*K_Y + (d-1)LP.$$ 

\par \vskip 1pc 
We now prove Theorem 2.1. Let $f, X$ be as in Theorem 2.1. 
Let $E_{\ell}$ ($1 \le \ell \le m$) be all of the negative curves on $X$,
each of which is either a $(-1)$ or $(-2)$-curve (Lemma 1.4). 
Since $K_X^2 \le 7$, we have $m \ge 1$. 
By the proof of Lemma 1.2, we have $deg(f) = d^2$ and $f^*E = dE$ for 
every negative curve $E$. 
Write 
$$K_X = f^*K_X + (d-1) \sum_{\ell = 1}^m E_{\ell} + Z,$$ 
where $Z$ is effective and supported by curves of non-negative self
intersection. 
Set $D = -K_X - \sum E_{\ell}$. Write $Z = \sum_i z_i Z_i$ with 
$z_i \ge 1$ and $Z_i$ irreducible. We first prove: 

\par \vskip 1pc 
{\bf Lemma 2.2} 

\par \noindent 
(1) $f^*C \sim dC$ for every divisor $C$ on $X$. 

\par \noindent 
(2) $-K_X = \sum_{\ell=1}^m E_{\ell} + D$. 

\par \noindent 
(3) {\it $Z \sim (d-1)D$ and $D$ is nef. 
We have $h^0(X, D) = 1 + (D. D - K_X)/2 \ge 1$,
so we will assume $D \ge 0$ thereafter. If $D \ne 0$,
then $h^0(X, D) = D^2 + 1 + e \ge 2$, 
where $e$ is the number of connected components of $\sum E_{\ell}$.} 

\par \noindent 
(4) {\it If $D \ne 0$ (i.e., if $Z \ne 0$) 
but $D$ is not big then $Z_i \cong {\bold P}^1$ with $Z_i^2 = 0$ 
(so $Z_i \sim Z_1$) and $D \sim e Z_1$.} 

\par \noindent 
(5) {\it If $D$ is big, then $|D|$ is base point free and
$dim |D| = D^2 + 1$.}

\par \vskip 1pc 
{\it Proof of Lemma $2.2$.} Since $K_X^2 \le 7$, it is easy 
to see that $Pic X$ is generated by negative curves $E$; 
now (1) follows from the fact that $f^*E = dE$. (2) is just the definition
of $D$. 
The first part of (3) follows from (1) and the ramification divisor formula
displayed above. 
Note that $Z_i^2 \ge 0$ as mentioned early. 
So $Z$ and hence $D$ are nef. By the Kawamata-Viehweg vanishing [K, V] we
have 
$H^1(X, D) = 0$ since $D - K_X$ is nef and big. 
Clearly, we have $H^2(X, D) \cong H^0(X, K_X-D) = 0$ since
$(d-1) D \sim Z \ge 0$ and the Kodaira dimension $\kappa(X)$
is negative.
So by the Riemann-Roch theorem, $h^0(X, D) = 1 + (D. D-K_X)/2 \ge 1$. 

\par
Assume $D \ne 0$. Then
$$\gather
D . \sum E_{\ell} = - \sum_{\ell} \{E_{\ell} . (K_X + \sum_p E_p)\} =
\sum_{\ell} (2 - \sum_{p \ne \ell} E_{\ell} . E_p) = \\
2m - \sum_{\ell \ne p} E_{\ell} . E_p = 2m - 2(m - e) = 2e. 
\endgather$$ 
In the second last equality, we use the fact that $\sum E_{\ell}$ 
is a union of $e$ connected rational trees for $|K_X + \sum E_{\ell}| = |-D|
= \emptyset$ 
and by Lemma 1.3, and that the number of edges in a connected tree is one
less 
than its number of vertices. 
Now $-D . K_X = D^2 + D . \sum E_{\ell}$. 
So $h^0(X, D) = 1 + D^2 + e$. This proves (3).

\par 
Suppose that $D \ne 0$ but $D$ is not big. So $Z^2 = 0 = D^2$. Thus
$Z_{\text{\rm red}}$ 
is a disjoint union of curves $Z_i$ with $Z_i^2 = 0$. Since $-K_X$
is nef and big, 
we have $K_X . Z_i < 0$ (equality would imply that $Z_i^2 < 0$ by the Hodge 
index theorem). So $Z_i \cong {\bold P}^1$ and $K_X . Z_i = -2$ by the genus
formula for $Z_i$. Now for all $i \ne j$, $Z_i \sim Z_j$ because $Z_i \cap
Z_j = \emptyset$. 
Since $D . Z = (d-1)D^2 = 0$, 
$D$ is contained in fibres of the ${\bold P}^1$-fibration
${\Phi}_{|Z_1|} : X \rightarrow {\bold P}^1$. 
This and $D^2 = 0$ imply that $D$ is a positive multiple of $Z_1$. 
Since the fibration clearly has a cross-section, 
$D \sim s Z_1$ for some positive integer $s$. 
Now $dim |D| = s$. This and (3) imply that $s = e$. This completes
the proof of (4).

\par 
Suppose that $D$ is big (and also nef). Then $dim |D| \ge 2$ by (3). 
Write $|D| = |M| + F$, where $F$ is the full fixed part. 
The Stein factorization applied to our regular surface $X$ implies that 
$M = kM_1$ where $M_1$ is irreducible. 
Since $|K_X+M_1| \le |K_X+D| = \emptyset$, the $M_1 \cong {\bold P}^1$ 
by Lemma 1.3. Thus $|M_1|$ is base point free and $dim |M_1| = M_1^2 + 1$ 
(see e.g. [DZ, Lemma 1.7]). Hence we may assume that $k = 1$ unless $M^2 = 0$. 
So either $M^2 = 0$ and $dim|M| = dim|kM_1| = k$, or $M^2 > 0$, 
$k = 1$ and $dim |M| = M^2 + 1$. 

\par 
Note that $D^2 = k^2M_1^2 + kM_1. F + D . F$, where each term of the latter 
is non-negative. 
If $M^2 = 0$, then the big and nefness of $D$ implies that $F > 0$ 
and ($D$ is 1-connected and hence) $M . F > 0$, whence $D^2 \ge k M_1 . F
\ge k$. 
This leads to $k = dim |M| = dim |D| = D^2 + e \ge k + 1$, a contradiction. 
Thus $M^2 > 0$ and hence 
$M^2 + 1 = dim |M| = dim |D| = D^2 + e \ge M^2 + M . F + 1$, 
whence $F = 0$ and $e = 1$. This proves (5). 

\par \vskip 1pc 
{\bf Lemma 2.3.} 

\par \noindent 
(1) {\it $D + \sum_{\ell=1}^m E_{\ell}$ is a simple loop of ${\bold P}^1$'s ($D$
being general), 
is a member of $|-K_X|$ and contains all negative curves on $X$.}

\par \noindent
(2) {\it $\sum E_{\ell}$ is connected with $m \ge 3$ components, i.e., $e = 1$.}

\par \noindent 
(3) {\it $LP := D + \sum E_{\ell}$ fits one of the $12$ Figures attached
at the end of the paper;
when $D^2 > 0$, $D$ is replaced by a member $D_1+D_2$ in $|D|$ with
$D_j^2 \ge 0$, $D_1 . D_2 = 1$ and simple normal crossing $LP := D_1 + D_2 +
\sum E_{\ell}$.}

\par \noindent
(4) {\it There is a smooth blow down $\mu : X \rightarrow {\bold P}^2$ of
curves
in $\sum E_{\ell} \le LP$ such that $\mu(LP)$ is a union of three
non-concurrent
lines $\sum L_j$ and $\mu$ is the blow up of some of the three intersection
points
of $\sum L_j$ and their infinitely near points. Moreover, the self-map $f$
of
$X$ induces a self-map of ${\bold P}^2$, also denoted by $f$ such that
$f \mu = \mu f$.}

\par \vskip 1pc 
{\it Proof of Lemma $2.3$.} (1) Since $K_X^2 \le 7$ and by Lemma 1.4, 
there is a birational morphism
$X \rightarrow {\bold P}^2$ (see, e.g. [DZ, Lemma 4.2]), and also
$m = \#(\sum E_{\ell}) \ge 3$; it is also easy to check that
if $m = 3$ then $K_X^2 = 7$ and $\sum E_{\ell}$ is a linear chain.
After choosing $D$ general (if $D \ne 0$), we may
assume that $D + \sum E_{\ell}$ ($\sim -K_X$) is reduced and of simple normal
crossing, and by Lemma 1.3 it is a simple loop of ${\bold P}^1$'s. So (1) is true.

\par
For (2), it is clear when $D = 0$ or $D$ is big 
(see the proof of Lemma 2.2).
Suppose that $D \ne 0$ but $D^2 = 0$. In notation of Lemma 2.2,
we have $Z_1 . \sum E_{\ell} = Z_1 . (-K_X - D) = 2$, say 
$Z_1 . \sum_{j=1}^k E_j = 2$ ($k = 1, 2$). So if $e \ge 2$ then $Z(1) +
Z(2) + \sum_{j=1}^k E_j$ contains a loop and is a sub-divisor of
the simple loop  $Z_1 + \cdots + Z(e) + \sum_{\ell = 1}^m E_{\ell}$, 
whence $m = k \le 2$,
a contradiction; here
$\sum_{j=1}^e  Z(j) \sim eZ_1$ is a general members of $D$.
So $e = 1$.

\par
To prove (3), we consider the case $D = 0$, or $D \ne 0$ and $D^2 = 0$,
or $D^2 > 0$ separately.
If $D \ne 0$ but $D^2 = 0$, we let $\varphi = \Phi_{|D|}$;
if $D = 0$ we let $\varphi : X \rightarrow {\bold P}^1$ be a 
suitable ${\bold P}^1$-fibration.
Let $S_1, \dots, S_k$ be all of the singular fibres [Z1, Lemma 1.5]. Clearly,
each $S_j$ consists of negative curves, so $\text{\rm Supp} \sum_j S_j \le \sum
E_{\ell}$.
We may assume that $E_1$ and $E_2$ in $\sum_{\ell=1}^m E_{\ell}$ are two
distinct cross-sections of $\varphi$; this is true when $D \ne 0$ but $D^2 = 0$,
and can be done by choosing
$\varphi$ suitably and making use of (1) and the fact that $m \ge 3$.

\par
Suppose that $D = 0$.
If $k \ge 2$, (1) implies that $k = 2$ and $LP = \sum E_{\ell}$ 
fits one of Figures 1 - 5.
If $k = 1$, then (1) implies that $E_1 . E_2 = 1$ and
$S_1$ is a linear chain of length
$3 - E_1^2 - E_2^2$ and dual graph
$$(-1) -- (-2) -- \cdots -- (-2) -- (-1).$$
Now the consideration of the new fibration $\varphi_2$ with a fibre $T_1$
formed by $E_1$ and $-E_1^2$ components in $S_1$, shows that $X$ has a new
negative curve in a fibre $T_2$ ($\ne T_1$) but outside $\sum E_{\ell}$, a contradiction.

\par
Suppose that $D \ne 0$ but $D^2 = 0$. 
Then (1) implies that $\varphi$ has exactly one singular fibre $S_1$,
and $LP := D + \sum E_{\ell}$ fits Figure 6, 7 or 8.

\par 
Suppose that $D^2 > 0$. Let $\mu : X \rightarrow {\bold P}^2$ be a
birational morphism [DZ, Lemma 4.2], smoothly blowing down curves in
$\sum E_{\ell} \le D + \sum E_{\ell} \sim -K_X$ so that 
$\mu_*(D + \sum E_{\ell})$ is a nodal member of $|-K_{{\bold P}^2}|$.
Then (1) implies that this nodal member is a union of a line
($= \mu(\sum E_{\ell}) = \mu(E_1)$ say) and a conic $\mu(D)$ so that $D +
\sum E_{\ell}$ fits Figure $j$ for some $9 \le j \le 12$,
where $D_1 + D_2$ should be read as $D$. 
Now let $D_1+D_2$ be the proper transform of the two lines meeting $\mu(E_1)$
at the two points $\mu(E_1) \cap \mu(D)$. Then from Figures 9-12 describing
$\mu$, we see that $D_1+D_2$ is a member in $|D|$ so that
$LP := D_1 + D_2 + \sum E_{\ell}$ is of simple normal crossing
and fits exactly Figure $j$. This proves (3).
(4) is clear from Figures 1 - 12 describing $\mu$, and Lemma 1.2.

\par \vskip 1pc
{\bf Lemma 2.4.} {\it For every $1 \le k \le 12$, there is exactly one
surface $X(k)$,
modulo isomorphism, with a simple loop $LP$ of Figure $k$ and satisfying
Lemma $2.3 (1)$, i.e., satisfying $(2a) - (2c)$ of Theorem $2.1$.}

\par \vskip 1pc
{\it Proof of Lemma $2.4$.} Let $X(j)$ denote a surface with $LP \sim -K_{X(j)}$
of Figure $j$. We see easily that there are following birational morphisms:
$$\gather
X(5) \rightarrow X(4) \rightarrow X(2) \rightarrow X(1) \leftarrow X(2) \leftarrow X(3), \\
X(6) \leftarrow X(7) \leftarrow X(8), \\
X(9) \leftarrow X(10) \rightarrow X(11) \leftarrow X(12).
\endgather$$

\par
We will prove for some cases with smaller $K_X^2$ and the others
either follow or are easier.
For $X(5)$, we can take $\mu : X(5) \rightarrow {\bold P}^2$ to be the
composite of
the blow up of the three intersection points $P_k$ ($k = 1, 2, 3$) of the
three non-concurrent lines $\sum_{i<j} L_{ij}$ (with $L_{ij}$ the line
joining $P_i$ and $P_j$, and
we may assume that $P_1 = [1:0:0], P_2 = [0:1:0], P_3 = [0:0:1]$)
with $E_k$ the exceptional divisor, and the blow up of the three
intersection
points of $E_1, E_2, E_3$ with respectively the proper transforms of
$L_{12}, L_{23}$ and $L_{31}$.

\par
For $X(8)$, we see that there is a birational morphism $\eta : X(8)
\rightarrow {\bold P}^1 \times {\bold P}^1$. It is the composite of the blow
up of two points $P_j$
having coordinates $(x, y) = (0, \infty)$ and $(\infty, \infty)$ with
exceptional divisors $E_j$, and the blow up of
the two intersection points of $E_1, E_2$ with respectively the proper
transforms
of the fibres $x = 0$ and $x = \infty$.

\par
For $X(10)$, we can take $\mu : X(10) \rightarrow {\bold P}^2$ to be the
composite of
the blow up of the two points $P_j$ ($j = 1, 2$) with exceptional divisor
$E_j$,
and the blow up of the intersection point of $E_1$ with the proper transform
of $L_{12}$
(we use the notation of $X(5)$). 
For $X(12)$, we can take $\mu : X \rightarrow {\bold P}^2$ to be the composite
of the blow up of $P_1$ with $E_1$ the exceptional curve, the blow up of the
intersection
point of $E_1$ with the proper transform of $L_{12}$ with $F_1$ the
exceptional divisor, and the blow up of the intersection point of $F_1$ with
the proper transform
of $L_{12}$ (we use the notation of $X(5)$).

\par
The above arguments shows that for each of $k = 5, 8, 10, 12$,
the surface $X(k)$ is unique modulo isomorphism. The other cases are
similar. This proves Lemma 2.4.

\par \vskip 1pc
Now we will complete the proof of Theorem 2.1. In view of Lemma 2.3,
we only need to prove Theorem 2.1 (4) and (5). Theorem 2.1 (4)
follows from Lemma 2.4. Indeed, if $\#\sum P_j = 0$ (resp. $1$, or $2$),
then as in Lemma 2.3, one can show that $LP$ which satisfies (2a) - (2c), 
fits one of Figures 1-4 (resp. Figures 6-8, or Figures 9-12).

\par
For Theorem 2.1 (5), we note that $-K_Y \sim LP$, and
the latter is given in one of Figures 1 - 12 and hence is nef and
big
so that $K_Y^2 = (LP)^2$ is between $3$ and $7$. In particular,
$Y$ is a weak del Pezzo surface. The last assertion in Theorem 2.1 (5) now 
follows from Example 1.5 with $g := \varphi_d$.
This proves Theorem 2.1 and also Theorem 2. Now Theorems 1 and 3
are consequences of Theorem 2.1 and Lemma 1.1.

\par \vskip 2pc
$$\gather
(-1) -- (-1) -- (-1) \\
\hskip 0.1pc |  \hskip 6.2pc  | \\
\hskip 0.1pc |  \hskip 6.2pc  | \\
(-1) -- (-1) -- (-1) \\
\text{\rm Figure 1}
\endgather $$

$$\gather
(-2) ---- (-1) ---- (-2) \\
\hskip 0.1pc |  \hskip 10pc  | \\
\hskip 0.1pc |  \hskip 10pc  | \\
(-1) -- (-1) -- (-1)-- (-1) \\
\text{\rm Figure 2}
\endgather $$

$$\gather
(-2) -- (-1) -- (-2) -- (-1) \\
\hskip 0.1pc |  \hskip 10pc  | \\
\hskip 0.1pc |  \hskip 10pc  | \\
(-1) -- (-2) -- (-1)-- (-2) \\
\text{\rm Figure 3}
\endgather $$

$$\gather
(-1) -- (-2) -- (-2) -- (-1) \\
\hskip 0.1pc |  \hskip 10pc  | \\
\hskip 0.1pc |  \hskip 10pc  | \\
(-2) -- (-1) -- (-1)-- (-2) \\
\text{\rm Figure 4}
\endgather $$

$$\gather
(-1) --- (-2) --- (-2) --- (-1) \\
\hskip 0.1pc |  \hskip 13pc  | \\
\hskip 0.1pc |  \hskip 13pc  | \\
(-2) -- (-2) -- (-1)-- (-2) -- (-2) \\
\text{\rm Figure 5}
\endgather $$

$$\gather
\hskip 0.1pc D \\
(-1) -- (0) -- (-1) \\
\hskip 0.1pc |  \hskip 6.2pc  | \\
\hskip 0.1pc |  \hskip 6.2pc  | \\
(-1) -- (-2) -- (-1) \\
\text{\rm Figure 6}
\endgather $$

$$\gather
\hskip 0.1pc D \\
(-1) ---- (0) ---- (-2) \\
\hskip 0.1pc |  \hskip 10pc  | \\
\hskip 0.1pc |  \hskip 10pc  | \\
(-1) -- (-2) -- (-2)-- (-1) \\
\text{\rm Figure 7}
\endgather $$

$$\gather
\hskip 2.7pc D \hskip 6pc \\
(-2) -- (0) -- (-2) -- (-1) \\
\hskip 0.1pc |  \hskip 10pc  | \\
\hskip 0.1pc |  \hskip 10pc  | \\
(-1) -- (-2) -- (-2)-- (-2) \\
\text{\rm Figure 8}
\endgather $$

$$\gather
\hskip 2.5pc D_1 \hskip 5pc D_2 \hskip 2.5pc \\
(0) ------ (0) \\
\hskip 0.1pc |  \hskip 6.2pc  | \\
\hskip 0.1pc |  \hskip 6.2pc  | \\
(-1) -- (-1) -- (-1) \\
\text{\rm Figure 9}
\endgather $$

$$\gather
\hskip 5pc D_1 \hskip 1.5pc D_2 \hskip 1.8pc \\
(-2) -- (0) -- (0) \\
\hskip 0.1pc |  \hskip 6.2pc  | \\
\hskip 0.1pc |  \hskip 6.2pc  | \\
(-1) -- (-2) -- (-1) \\
\text{\rm Figure 10}
\endgather $$

$$\gather
\hskip 2.5pc D_1 \hskip 5pc D_2 \hskip 2.5pc \\
(0) ------ (1) \\
\hskip 0.1pc |  \hskip 6.2pc  | \\
\hskip 0.1pc |  \hskip 6.2pc  | \\
(-2) -- (-1) -- (-1) \\
\text{\rm Figure 11}
\endgather $$

$$\gather
\hskip 5pc D_1 \hskip 1.5pc D_2 \hskip 1.8pc \\
(-2) -- (0) -- (1) \\
\hskip 0.1pc |  \hskip 6.2pc  | \\
\hskip 0.1pc |  \hskip 6.2pc  | \\
(-2) -- (-1) -- (-2) \\
\text{\rm Figure 12}
\endgather $$

\head
References
\endhead

\par \noindent 
[B] A. Beauville, Endomorphisms of hypersurfaces and other manifolds,
math.AG / 0008205. 

\par \noindent
[Da] M. Dabija, Algebraic and geometric dynamics in several complex variables,
Ph D thesis, Univ. of Michigan, 2000.

\par \noindent 
[D] I. Dolgachev, Personal e-mail dated 24 December 2001.

\par \noindent
[DZ] I. Dolgachev and D. -Q. Zhang, Coble rational surfaces,
Amer. J. Math. {\bf 123} (2001), no. 1, 79--114.

\par \noindent 
[GZ1, 2] R. V. Gurjar and D. -Q. Zhang, 
$\pi_1$ of smooth points of a log del Pezzo surface is
finite, I; II, J. Math. Sci. Univ. Tokyo {\bf 1} (1994), no. 1, 137--180;
J. Math. Sci. Univ. Tokyo {\bf 2} (1995), no. 1, 165--196.

\par \noindent  
[GPZ] R. V. Gurjar, C. R. Pradeep and D. -Q. Zhang,
On Gorenstein Surfaces Dominated by ${\bold P}^2$,
Nagoya Mathematical Journal, to appear,
math.AG/0112242.

\par \noindent 
[K] Y. Kawamata, A generalization of Kodaira-Ramanujam's vanishing theorem,
Math. Ann. {\bf 261} (1982), no. 1, 43--46. 

\par \noindent 
[MZ1, 2] M. Miyanishi and D. -Q. Zhang,
Gorenstein log del Pezzo surfaces of rank one, I; II, J. Algebra
{\bf 118} (1988), no. 1, 63--84; J. Algebra {\bf 156} (1993), no. 1, 183--193. 

\par \noindent 
[M] D. Mumford, 
The topology of normal singularities of an algebraic surface and a
criterion for simplicity, Inst. Hautes Itudes Sci, Publ. Math. 
No. {\bf 9} (1961), 5--22. 

\par \noindent 
[R] M. Reid, Chapters on algebraic surfaces, in:
Complex algebraic geometry (Park City, UT, 1993), 3--159, 
IAS/Park City Math. Ser., {\bf 3}.

\par \noindent 
[V] E. Viehweg, Vanishing theorems, J. Reine Angew. Math. {\bf 335} (1982), 1--8.

\par \noindent 
[Y] Q. Ye, On Gorenstein log del Pezzo Surfaces,     
Japanese J. Math. to appear, math.AG / 0109223.

\par \noindent 
[Z1] D. -Q. Zhang, Logarithmic del Pezzo surfaces of rank one 
with contractible boundaries, Osaka J. Math. {\bf 25} (1988), no. 2, 461--497.

\par \noindent 
[Z2] D. -Q. Zhang, Logarithmic del Pezzo surfaces with rational double 
and triple singular points, Tohoku Math. J. {\bf 41} (1989), no. 3, 399--452.

\par \noindent 
[Z3] D. -Q. Zhang, Algebraic surfaces with nef and big anti-canonical divisor, 
Math. Proc. Cambridge Philos. Soc. {\bf 117} (1995), no. 1, 161--163. 

\par \noindent
[ZD] D. -Q. Zhang, Automorphisms of finite order on rational surfaces,
with an appendix by I. Dolgachev, J. Alg. {\bf 238} (2001), 560--589.

\par \noindent
[N] N. Nakayama, Ruled surfaces with non-trivial surjective endomorphisms,
RIMS, Kyoto Univ. Preprint No. 1286, July 2001, Kyushu J. Math. to appear.

\par \noindent
[S] M. Segami, On surjective endomorphism of surfaces (in Japanese),
Proceedings of Symposium on vector bundles and algebraic geometry,
Jan. 1997, Kyushu Univ. ed. S. Mukai and E. Sato, 93-102.

\par \noindent
[F] Y. Fujimoto, Endomorphisms of smooth projective threefolds with non-negative
Kodaira dimension, Preprint 2000.

\par \noindent
[FS] Y. Fujimoto and E. Sato, On smooth projective threefolds with non-trivial
surjective endomorphisms, Proc. Japan Acad. {\bf 74} (1998), 143 -- 145.

\enddocument